# Romancing Mathematics with Chemistry—How Mathematical Trees Can Be Used to Synthesize Molecular Structures

**Chin-yah Yeh**
*Utah Valley University*

## Abstract

*Structures of chemical compounds can be synthesized and categorized through mathematical means. Organic compounds are suitable targets because of their simple valences. Acyclic organic compounds made of hydrogen and second-row elements C, N, O, and F are presented as an example. In five categories of organic compounds, chemical structures can be generated exclusively and exhaustively using ab initio methods. It is shown that mathematical variables can serve as chemical symbols and mathematical equations are chemical structure generators.*

# I. BACKGROUND

Chemists have a long tradition of using atomic valences to find molecular structures graphically. Chemistry relies on this simple procedure to make progress; however, drawing does not generate structures systematically, i.e., there is no guaranteed completeness or uniqueness. Are all structures generated? Is a structure generated more than once? Take butane molecules ($C_4H_{10}$) for instance. It can be seen (Fig. 1) that one straight-chain and one branched structure suffice. But as molecular size gets bigger, finding structures exhaustively and uniquely can become a daunting job. This paper deals with how we generate chemical structures from first principles. Compared with the empirical method that relies on computers,[1] structural generation through first principle has a long and winding past.

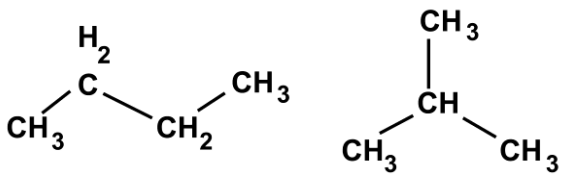

**Fig. 1.** Structures for butanes

About the time Darwin published his *Origin of Species*, two mathematicians endeavored in structuring chemical compounds. They are called the 'invariant twins,'[2] Sylvester and Cayley. Their idea was to start from two observations. First, certain algebraic equations do not vary under geometrical transformations (a circle does not vary under rotation around an axis through its center, for instance). Secondly, chemical structures and molecular properties do not vary under the same transformation either. What is the connection? Sylvester related chemical structures, as understood at his time, to invariant algebraic forms. The result was published[3] in 1878 but does not lead very far. Cayley[4] acknowledged Sylvester's effort and focused on simple acyclic structures. Cayley's invention of an analytical form called trees has proven to be the right tool for generating chemical structures, even though he used this analysis only for counting the structures.

# II. CRITICISM AND POTENTIAL OF CAYLEY'S METHOD

Cayley's result[5] caused numerous revisions shortly afterwards,[6,7] partly because of some errors found in his paper, as he used a manageable and yet laborious method that separates trees (which consist of lines and nodes) into centric and bicentric ones (Fig. 2). Nonetheless, what Cayley established should be seen as more of an analytical tool than a calculation. Some 70 years later, his method was proven rigorously by Otter[8] and generalized by Harary and Norman[9] (see Section IV). These authors improved Cayley's method by bypassing centric and bicentric trees. As a result, the method becomes harder to carry out manually, but easier to code (because of more repetitive steps; see Yeh[10] on coding). The present paper points out that Cayley and the above-mentioned authors have unknowingly (or unreportedly) discovered a mathematical tool for synthesizing a good part of the 10–20 million chemical compounds known today.

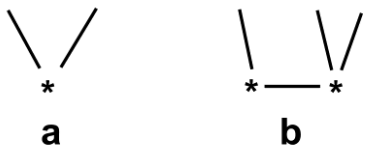

**Fig. 2.** Examples of (a) centric trees and (b) bicentric trees.

Inspired by Sylvester's work on the interchange of variables in differential calculus, Cayley[4] started the concept of mathematical trees. In modern terms, rooted trees are enumerated by the series expansion of a function $f(x)$ of an independent variable $x$ as expressed in the equation

$$f(x) = xe^{f(x)}. \quad (1)$$

Then Otter's formula[8] $F(x) = f(x) - \frac{1}{2}\{[f(x)]^2 - f(x^2)\}$ is used to extract (or unlabel) root-free trees (see Fig. 3). Cayley did not limit the number of lines connecting to a node in a tree to four, but four happens to be the valence of a carbon atom in all organic compounds. We are able to adopt Cayley's scheme in organic compounds and include other essential atoms such as nitrogen, oxygen, and halogens. In Sections III and IV, we shall use Cayley's scheme on the simplest class of organic molecules, the alkane series ($C_nH_{2n+2}$), to illustrate two points. One is to show how his scheme complies with the concept of chemical structure, although he himself treated the alkane series differently—by counting centric and bicentric trees. The other is the central issue of this paper: Chemical elements can serve as mathematical variables and equations are generators of chemical formulas.

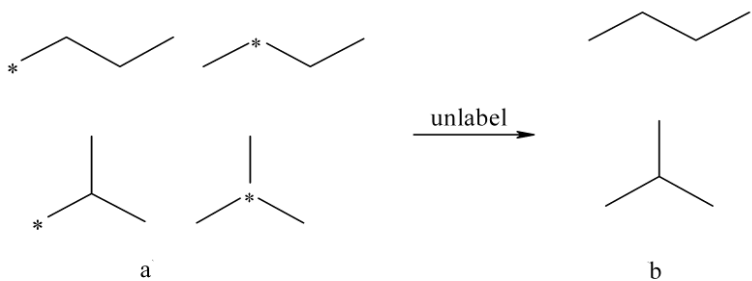

**Fig. 3.** Examples of roots in trees: (a) rooted; (b) root-free.

In essence, our instinct to draw does not quite work on chemical structures. The right strategy is to take a step back and find all rooted (labeled) structures first, as shown in Section III, and then unlabel them, as shown in Section IV.

## III. ROOTED ALKANE SERIES ($C_nH_{2n+2}$)

Nothing is new under the sun. Today there are 10–20 million known chemical compounds, 90% of which are organic. A good part of these organic compounds were already programmed in the 19$^{th}$ century, at least mathematically, by the analytical forms called trees.[11]

Trees are basic constituents of chemical structures. Cayley used trees and the usual algebra to count acyclic molecules; however, trees differ from most algebras in that operations are primitive, i.e., they have fewer properties than most algebras. Ordinary numbers are immune to change of order or association, but when these operations are applied to elements of trees, a distinct tree results. The domain of trees is ideal to represent chemical elements with which the molecules grow. We harvest the technique developed for tree enumeration in number operations and use it on the growth of trees. That is, the equation for counting, eq. (1), is also used for tree generation. This grafting works exceedingly well. When $x$ is a number, trees of the same size (same number of nodes) are lumped together, but when $x$ is a node of a tree, each tree is uniquely represented. Variable $x$ is a chemical element and thus has lost all versatile features of a number.

Chemical structures are more restricted than mathematical trees in two aspects. First, nodes in chemical formulas are structured, whereas those in trees are not. Carbon atoms are tetravalent and tetrahedral in shape for saturated compounds (no multiple bonds). Second, if a carbon atom is bonded to less than four other carbons, it is padded with hydrogens. With the first aspect taken into account, eq. (1) can be reduced to[12]

$$a(x) = 1 + (x/6)\{[a(x)]^3 + 3a(x)a(x^2) + 2a(x^3)\}. \quad (2)$$

How does eq. (2) become a generator of chemical formulas? We view the growth of trees from two angles, one from the equation and the other from its solution $a(x)$. Eq. (2) reveals the growth process, whereas its solution gives individual tree structures. Key for solving eq. (2) is iteration. Meanwhile, a homologous series (such as the alkane series) of acyclic chemical compounds is generated graphically through iteration. Structures are produced exclusively and exhaustively by first principle. Iteration is carried out by approximating $a(x)$ successively as $a_0$, $a_1$, $a_2$, etc. As in trees, when $x$ is a number, eq. (2) is an ordinary algebraic equation used for enumeration. When $x$ is an atom, $a(x)$ represents a series of chemical structures. Hydrogen padding is carried out by setting $a_0$ equal to H. The first term in eq. (2) namely 1, is symbolized as a hydrogen atom H and used as the iteration initiator $a_0$; each $x$ is a carbon atom C. Now $a(x)$ becomes a function of two variables

$$a(H,C) = H + (C/6)\{[a(H,C)]^3 + 3a(H,C)a(H^2,C^2) + 2a(H^3,C^3)\},$$

and is iterated as follows,

$$a_{n+1}(H,C) = H + (C/6)\{[a_n(H,C)]^3 + 3a_n(H,C)a_n(H^2,C^2) + 2a_n(H^3,C^3)\},$$

leading to the solution

$$a_0 = \cdot H$$
$$a_1 = \cdot H + \cdot CH^3$$
$$a_2 = \cdot H + (\cdot C/6)[(H + CH^3)^3 + 3(H + CH^3)(H^2 + (CH^3)^2) + 2(H^3 + (CH^3)^3)]$$
$$\quad = \cdot H + \cdot CH^3 + \cdot CH^2CH^3 + \text{higher terms}$$
$$a_3 = \cdot H + \cdot CH^3 + \cdot CH^2CH^3 + \{\cdot CH^2CH^2CH^3 + \cdot CH(CH^3)^2\} + \text{higher terms}$$

where dots · are added in front of the roots H or C to represent the radical nature of these atoms. Sequence $a_n$ grows at each successive iteration by adding terms of $n$th order, which come from linking lower-order terms to the root. Merely by lowering superscripts into subscripts, sequence $a_n$ becomes alkyl radicals ($C_nH_{2n+1}$) collectively. Iteration also modifies the meaning of the equal sign, which now stands for 'replaced by,' as used in programming language.

Four remarks on eq. (2) are in order. First, iteration is a growth process, but molecules viewed as mathematical trees do not grow exactly as biological trees. Existing mathematical trees are used as branches and bonded to the root to form new trees. Molecules grow only from the root (not from other nodes) and do not grow bigger; they grow by forming new alkane molecules. Each carbon atom has four valences, with one pointing to the root and three pointing away; we can call them one stem and three branches. The root has three branches and no stem. Second, as variables in eq. (2) are reinterpreted as atoms C and H, operations become primitive. This results in terms that are all distinct from one another, exactly what a proper structure generator needs. For instance, there are two distinct $C^3H^7$ terms, representing two propyl radicals. Third, each carbon atom is tetrahedral and hence has three equivalent sites for bonding. Eq. (2) must obey the Pólya enumeration theorem[13] and engender molecular group symmetry $C_{3v}$. Note that the symmetry applies to the sites but not necessarily to the actual molecules. The three terms in the curly brackets of eq. (2) correspond to elements of $C_{3v}$, with coefficients 1, 3, and 2, the groupings in $C_{3v}$. Within the constraint of tree operations and node symmetry $C_{3v}$, eq. (2) generates unique structures. Nodes of nitrogen or oxygen atoms are even simpler and can be appended to eq. (2) as seen in Section V. Fourth, symmetry reduces the number of molecular structures. Pólya's theorem rationalizes and quantifies this fact.

## IV. FREE ALKANE SERIES ($C_nH_{2n+2}$)

In each structure generated by eq. (2), the root is affiliated with the radical quality of each alkyl radical. Our next task is to remove the root. For trees with no constraint from chemical valence, the roots are removed by using Otter's formula.[8] But chemical structures are either a subset of trees, such as the alkane series, or its extended version, such as the organic compounds shown in the next section, with more than one kind of root for growth. For these cases, or any other case except that of genuine mathematical trees (made of structureless lines and nodes), Otter's formula is no longer right and has to be replaced by a more general form, that of the dissimilarity characteristic theorem[9] (DCT). DCT is not limited to graphs with the constraint of chemical valence. A big plus of DCT is that it survives converting the domain of variables from numbers to chemical elements. In other words, the expression of DCT is also a structure generator. Besides, DCT is a generalization of Euler's characterization theorem (see below in this section). The essential idea of DCT is this. Contrary to what Fig. 3 implies, root-free trees do not come from literally unlabeling rooted species, but from combining rooted trees of smaller size. Two rooted trees are combined to lose their roots, resembling two radicals/spins combined to form a non-radical with no spin. The result is a root-free (or free) molecule. Radicals outnumber free molecules, as seen in Fig. 3, because a molecule can be built from two radicals in many ways. When two radicals are combined, radicality (the label * in the graph) could remain on either side of the link or, alternatively, on the link itself, resulting in two node-labeled molecules and a link-labeled one. Every node-labeled molecule cancels a link-labeled one numerically. After cancellation, the amount left is the same as the number of free molecules. If the number of node-labeled molecules is $p$ and link-labeled ones $q$, the count of free molecules $\varphi$ is $\varphi = p - q$. This result is correct only if the two parts are dissimilar from each other. If they are the same, there is only one node-labeled species. Therefore, the net count should be corrected by adding a link-centered species, namely, a molecule made of two equal halves. The above argument works for multi-link molecules as well as single-link ones (diatomics). The count of free molecules is the difference between node-labeled species and link-labeled species with a correction whenever link-centered species $r$ exist. Namely, the net count should be $\varphi = p - q + r$. Consider butanes for example (Fig. 4).

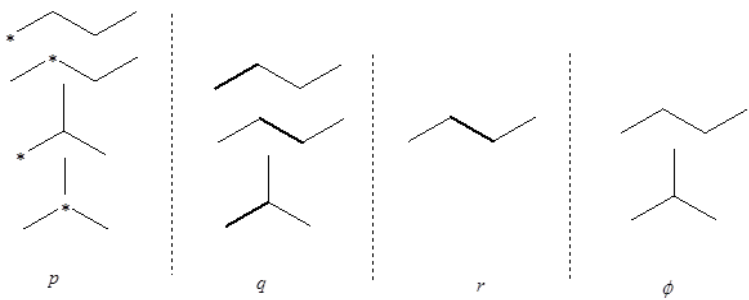

**Fig. 4.** Unlabeling butanes: $p$ (node-labeled) − $q$ (link-labeled) + $r$ (center-labeled) = $\varphi$ (label-free).

Let us imaginatively label one node with the isotope $^{14}C$ in each butane. There are four node-labeled butanes. A link can be labeled with isotope $^{14}C$ on both sides of the link as shown in Fig. 4. There are three link-labeled species and a link-centered one. So, there are $4 - 3 + 1 = 2$ butanes. The structures represented by $p$, $q$, and $r$ are conveniently calculated from rooted species $a$, which now is a function of $H$ and $C$, as follows:

$$p(H,C) = (C/24)\{[a(H,C)]^4 + 6[a(H,C)]^2 a(H^2,C^2) + 8a(H,C)a(H^3,C^3)$$
$$+ 3[a(H^2,C^2)]^2 + 6a(H^4,C^4)\}$$

$$q(H,C) = (1/2)\{[a(H,C) - H]^2 + [a(H^2,C^2) - H^2]\}$$

$$r(H,C) = a(H^2,C^2) - H^2$$

and

$$\varphi(H,C) = p(H,C) - q(H,C) + r(H,C).$$

Thus, with the help of DCT, we have turned a mathematical equation like eq. (2) into a machine that produces chemical structures. Quite remarkably, eq. (2) produces structures uniquely and exhaustively.

Three remarks are in order. First, knowledge may not follow chronological order. Pólya's theorem and DCT were not available in Cayley's time. Nonetheless, Cayley progressed in the right direction and obtained the correct result (understandably with minor computational errors). DCT proves to be an essential part of the scheme if we want to generate chemical structures. Work[6] that shuns DCT and claims a correction over Cayley's result is done *ad hoc* for counting and not able to generate structures. Second, DCT is none other than an extension of Euler's theorem, which states that the number of node-labeled species $p$ is equal to that of link-labeled species $q$ plus one, for an acyclic molecule with irregular shape. But, in reality, a molecule may possess symmetry, which corrects the formula to $p - q + r = 1$. In the case of a class of molecules, $p - q + r$ gives a count of the molecules. For such calculation, a concise source code has been written.[10] Third, eq. (2) does not include chiral isomers. Chiral structures of alkane molecules ($C_nH_{2n+2}$) have been enumerated by a modified[12] eq. (2), which can become a structure generator, too.

We have used this scheme on five categories of organic compounds. They are

1. Saturated acyclic compounds made of hydrogen and four common second-row elements C, N, O, and F.
2. Saturated acyclic compounds made of hydrogen and second-row elements C, N, and O, with no (N/O)–(N/O) connection.
3. Saturated acyclic compounds made of hydrogen and second-row elements C, N, and O, with no geminal (N/O) branches.
4. Acyclic compounds with multiple bonds between carbon atoms.
5. Aldehydes and ketones. They can be considered as a variation of Category II in the form of dehydrated geminal di-alcohols.

Formulas of enumerating these categories have thus been found for the first time and then converted to structure generators. In each category, a rooted species is generated first, and then DCT is used to calculate free molecules. A rooted species is a radical but corresponds to a mono-substituted molecule with a substituent, such as a nitrile group or a halogen. Therefore, as a byproduct, any mono-substituted series can be enumerated or generated by carrying out the first portion of the procedure. The romance of chemistry with math has a fruitful production. We now show the result of Category 1 explicitly as an example of the method.

## V. SATURATED ACYCLIC COMPOUNDS MADE OF HYDROGEN AND SECOND-ROW ELEMENTS C, N, O, AND F

Organic compounds contain four essential elements: H, C, N, and O. Using Cayley's scheme, we could include any one of the four second-row elements C, N, O, or F as root for growing acyclic organic molecules. In other words, we need to add three more terms in eq. (2) besides the one covering carbon. Thus, at the root, carbon can grow three branches, nitrogen two, oxygen one, and fluorine none. In turn, each branch can start with C, N, or O. These two alternating actions form an endless loop to produce all combinations of acyclic structures. The governing equation becomes

$$\begin{aligned} a(H,C,N,O,F) = &\cdot H + (\cdot C/6)\{[a(H,C,N,O,F)]^3 + 3a(H,C,N,O,F) \\ &\times a(H^2,C^2,N^2,O^2,F^2) + 2a(H^3,C^3,N^3,O^3,F^3)\} \\ &+ (\cdot N/2)\{[a(H,C,N,O,F)]^2 + a(H^2,C^2,N^2,O^2,F^2)\} \\ &+ \cdot O a(H,C,N,O,F) + \cdot F. \end{aligned}$$

Again, as variables turn into chemical symbols, the equation becomes a structure generator. The number of structures is enormously larger than that of the alkane series. Algebraically, fluorine behaves as hydrogen. Therefore, iteration of eq. (3) starts with two elements, $\cdot H + \cdot F$, instead of hydrogen alone. With the padding of hydrogens and fluorines, the result gives all isomers of fluorine compounds. As a special case where only the $\cdot F$ term is added to the right side of eq. (2), we arrive at all fluorinated hydrocarbons. Terms covering a single root, be it C, N, O, or F, are simply additive. Versatility of Cayley's scheme is eminent. With fluorine being less prevalent in chemistry and to avoid cluttering the formulas, we shall drop all fluorine derivatives by setting $F$ to 0 in what follows. The terms that come out the first iteration in solving eq. (3) are none other than the common radicals shown below

$$a(H,C,N,O) = \cdot H + \cdot CH_3 + \cdot NH_2 + \cdot OH + \cdot CH_2CH_3 + \cdot CH_2NH_2 + \cdot CH_2OH$$
$$+ \cdot NHCH_3 + \cdot NHNH_2 + \cdot NHOH + \cdot OCH_3 + \cdot ONH_2 + \cdot OOH + ...$$

By going through the root-removing process using DCT, we arrive at the assembly of free molecules of Category 1 as

$$\varphi(H,C,N,O) = CH_4 + NH_3 + OH_2 + CH_3CH_3 + CH_3NH_2 + CH_3OH$$
$$+ NH_2NH_2 + NH_2OH + HOOH + ...$$

As easily seen, $\varphi(H,C,N,O)$ covers most plain organic molecules. Chemical structures are neatly laid out through *ab initio* calculation. When variables are set as real numbers, $H = 1$, $C = N = O = x$, and $F = 0$, the structure generators are converted to enumerators $a(1,x,x,x)$ [from eq. (3)] and $\varphi(1,x,x,x)$. For saturated acyclic compounds made of C, H, N, and O atoms, counts of rooted and free chemical structures ordered by size (corresponding to the sum of C, N, and O atoms, as shown in the exponents) are

$$a(1,x,x,x) = 1 + 3x + 9x^2 + 39x^3 + 181x^4 + 921x^5 + 4920x^6 + 27408x^7 + 156948x^8 + 919361x^9 + 5480371x^{10} + ...$$

and

$$\varphi(1,x,x,x) = 3x + 6x^2 + 18x^3 + 65x^4 + 258x^5 + 1140x^6 + 5436x^7 + 27262x^8 + 142311x^9 + 766073x^{10} + ...,$$

respectively. Each count of $x^n$ corresponds to a distinct chemical structure of size *n*. With eq. (3) as another example, a mathematical equation is turned into a chemical structure generator. Details will be presented elsewhere.

## VI. CONCLUSION AND PROSPECT

In five categories of organic compounds, we are able to generate chemical structures exclusively and exhaustively using an *ab initio* method. It is demonstrated that mathematical variables can serve as chemical symbols and mathematical equations are chemical structure generators. The method starts as a counting tool for tree structures but finds a much better use as a chemical structure generator. This *ab initio* structure generating tool has a prospective use: It is more direct, more intuitive, and less error-prone than the currently popular algorithmic methods[1] used for organizing chemical compounds in a chemical data bank. One drawback of Cayley's method is its failing on polycyclic structures. Static cyclic structures are covered by Pólya's theorem, but growth on cyclic structures similar to tree growth hits a wall. Seeking a breakthrough in this topic has haunted chemists ever since, in spite of some limited success in polyhex systems.[14]

## Acknowledgments

The author is very grateful to Professor Jean-Francois Van Huele for his keen interest in this paper. His help on style and clarification throughout the whole paper is greatly appreciated.